\documentclass[a4paper,11pt,twoside]{article}
\RequirePackage[a4paper,head=1cm]{geometry}
\usepackage[utf8x]{inputenc}
\usepackage{textcomp}
\usepackage{amsmath,amsfonts,verbatim,afterpage,euscript,mathrsfs,amssymb,empheq}
\usepackage{amsthm}
\usepackage{dsfont}
\usepackage{hyperref}
\usepackage{pst-fill,pst-grad}
\usepackage{textcomp}
\usepackage{multicol}
\usepackage{amsmath}
\usepackage{amssymb}
\usepackage{hyperref}
\usepackage{dsfont}
\usepackage{soul}

\newtheorem{Theoreme}{Theorem}
\newtheorem{Lemme}{Lemma}[section]
\newtheorem{Corollaire}{Corollary}[section]

\newcommand{\mysection}{\setcounter{equation}{0} \section}

\title{New weighted pointwise inequalities and applications to generalized Sobolev estimates} 
\author{Diego Chamorro\footnote{Laboratoire de Math\'ematiques et Mod\'elisation d'Evry (LaMME) - UMR 8071. Universit\'e d'Evry Val d'Essonne, 23 Boulevard de France, 91037 Evry Cedex, France. email: \textit{diego.chamorro@univ-evry.fr}}, Anca-Nicoleta Marcoci\footnote{Department of Mathematics and Computer Science. Technical University of Civil Engineering, Bucharest, Bld. Lacul Tei, no. 124, sector 2. Romania. email: \textit{anca.marcoci@utcb.ro}}, Liviu-Gabriel Marcoci\footnote{Department of Mathematics and Computer Science. Technical University of Civil Engineering, Bucharest, Bld. Lacul Tei, no. 124, sector 2. Romania.  email: \textit{liviu.marcoci@utcb.ro}}.} 
\begin{document} 
\maketitle 
\begin{scriptsize}
\abstract{ \noindent In this article we study some new pointwise inequalities between rough singular integral operators, weighted maximal functions of the gradient and weighted Morrey spaces. These pointwise estimates will naturally lead us to a new class of weighted Sobolev-type inequalities.}\\

{\footnotesize
\noindent \textbf{Keywords: Singular integral operators; pointwise estimates; Sobolev inequalities.} \\
\noindent \textbf{MSC (2020) Primary: 42B20; Secondary: 42B25}
}
\end{scriptsize}
\mysection{Introduction and presentation of the results}
We want to study in this article some generalizations of the following pointwise inequality: for $n\geq 2$, if $f:\mathbb{R}^n\longrightarrow \mathbb{R}$ is a function such that $f\in \mathcal{C}^\infty_0(\mathbb{R}^n)$, then we have
\begin{equation}\label{Ineq_IntroRiesz}
|f(x)|\leq CI_{1}(|\nabla f|)(x),
\end{equation}
where the operator $I_{1}$ corresponds to the usual Riesz potential defined by the expression
$$I_{1}(f)(x)= \pi ^{-n/2}2^{-1}\frac{\Gamma ((n-1 )/2)}{\Gamma (1/2)}\int_{\mathbb{R}^n}\frac{f(y)}{|x-y|^{n-1}}dy.$$
Following \cite{Hoang2}, we will call the estimate (\ref{Ineq_IntroRiesz}) a \emph{subrepresentation formula}. From this pointwise estimate, it is possible to deduce in a very straightforward manner some important functional inequalities: indeed, since the Riesz potential $I_{1}$ is bounded from $L^p(\mathbb{R}^n)$ to $L^{q}(\mathbb{R}^n)$ with $1<p<n$ and $\frac{1}{q}=\frac{1}{p}-\frac{1}{n}$ (see \cite[Theorem 6.1.3]{Grafakos}), we can easily obtain from the previous pointwise control (\ref{Ineq_IntroRiesz}) the functional inequality
$$\|f\|_{L^{q}}\leq C\|I_{1}(|\nabla f|)\|_{L^{q}}\leq C\|\nabla f\|_{L^p},$$
which is the classical Sobolev inequality $\|f\|_{L^{q}}\leq C\|\nabla f\|_{L^p}$.\\

Some variants of this inequality were recently studied in \cite{ChMarcociMarcoci2}, \cite{Hoang}, \cite{Hoang1}, \cite{Hoang2}, \cite{Li} (see also the references therein) and the modifications given in these articles follow essentially two directions. The first one replaces the function $|f|$ in the left-hand side of (\ref{Ineq_IntroRiesz}) by $|T(f)|$ where $T$ is a generic singular integral operator. This modification implies a different treatment of the pointwise estimate (\ref{Ineq_IntroRiesz}) but it also has important consequences when considering functional inequalities as we can not ``simply'' obtain the boundedness of this kind of operators in the Lebesgue spaces (\emph{i.e.} in full generality we do not have the control $\|T(f)\|_{L^q}\leq \|f\|_{L^q}$) and thus another approach must be performed to deduce Sobolev-like estimates. The second type of modifications of the subrepresentation formula (\ref{Ineq_IntroRiesz}) above consists in replacing the Riesz potential $I_{1}$ in the right-hand side by a more general operator $\mathcal{F}$ and again the arguments must be different in order to deduce the wished pointwise inequalities.\\

In this article we will generalize the following generic subrepresentation formula (see Theorem \ref{Theo_Pointwise_Weighted} below):
$$|T(f)(x)|\leq C\mathcal{F}(|\nabla f|)(x),$$
and we will also establish some related functional inequalities (see Corollary \ref{Corollary2} below). In order to present our results, we need to be more precise about the operators $T$ and $\mathcal{F}$ in the previous expression: indeed, the operator $T$ will be a rough singular integral operator (whose properties will be governed by a kernel $\Omega$) and the operator $\mathcal{F}$ will be a weighted fractional version of the Riesz potentials (whose properties will naturally depend on the weight considered). We thus have the following definitions:

\begin{itemize}
\item[$\bullet$] {\bf The Operators $T_\Omega$ and $T^*_\Omega$}. 
For a locally integrable function $f:\mathbb{R}^n\longrightarrow \mathbb{R}$, we consider the rough singular integral operator $T_\Omega$ given by the formula
\begin{equation}\label{Def_Operator}
T_\Omega(f)(x)=p.v.\int_{\mathbb{R}^n}\frac{\Omega(y/ |y|)}{|y|^n}f(x-y)dy,
\end{equation}
where $\Omega:\mathbb{S}^{n-1}\longrightarrow \mathbb{R}$ is a function such that $\Omega \in L^1(\mathbb{S}^{n-1})$ and such that $\displaystyle{\int_{\mathbb{S}^{n-1}}\Omega \ d\sigma=0}$. Associated to this operator $T_\Omega$ , we can study the maximal rough singular integral operator $T^*_{\Omega}$ defined by the following expression
\begin{equation}\label{Def_OperatorTstar}
T^*_{\Omega}(f)(x)=\underset{t>0}{\sup}\left|\int_{\{|y|>t\}}\frac{\Omega(y/|y|)}{|y|^n}f(x-y)dy\right|.
\end{equation}
Of course, the relationship between these operators is quite straightforward since we have estimate 
\begin{equation}\label{EstimateEtoile}
|T_\Omega(f)(x)|\leq T^*_{\Omega}(f)(x),
\end{equation}
and for this reason, here we will focus rather on the operator $T^*_{\Omega}$ than on the operator $T_{\Omega}$.\\
\end{itemize}
As it can be expected, the properties of the operators $T_{\Omega}$ and $T^*_{\Omega}$ are closely related to the information available on the kernel function $\Omega$, and the study of the boundedness of such operators in connection with subrepresentation formulas of the form (\ref{Ineq_IntroRiesz}) is the main topic of this article.\\

\noindent Indeed, if we assume moreover that $\Omega \in L^{\infty}(\mathbb{S}^{n-1})$ then the following estimate was proven in \cite{Li}:
\begin{equation}\label{Ineq_IntroOper1}
|T_{\Omega}(f)(x)|\leq C\|\Omega\|_{L^{\infty}(\mathbb{S}^{n-1})}I_{1}(|\nabla f|)(x)\qquad \mbox{for }  f\in\mathcal{C}^\infty_0(\mathbb{R}^n).
\end{equation}
In \cite{Hoang}, the inequality (\ref{Ineq_IntroOper1}) was improved by considering the more general condition $\Omega \in L^{n,\infty}(\mathbb{S}^{n-1})$ (where the space $L^{n,\infty}$ is a Lorentz space) and the pointwise estimate is then 
\begin{equation}\label{Ineq_IntroOperLorentz}
|T_{\Omega}(f)(x)|\leq C\|\Omega\|_{L^{n,\infty}(\mathbb{S}^{n-1})}I_{1}(|\nabla f|)(x).
\end{equation}
Let us mention that in our recent work \cite{ChMarcociMarcoci2} we generalized these two previous estimates by considering the condition $\Omega \in L^{\rho}(\mathbb{S}^{n-1})$ for some $1<\rho<n$. Indeed, in this particular setting we have the space inclusions 
\begin{equation}\label{SpaceInclusions0}
L^{\infty}(\mathbb{S}^{n-1})\subset L^{n,\infty}(\mathbb{S}^{n-1})\subset L^{\rho}(\mathbb{S}^{n-1}), 
\end{equation}
see Section \ref{Secc_FuncSpaces} for more details about these spaces inclusions.\\

In the work \cite{Hoang1}, the pointwise estimates (\ref{Ineq_IntroOper1}) and (\ref{Ineq_IntroOperLorentz}) were studied by introducing a Riesz potential of order $0<\alpha <n$. Indeed, for $0<\alpha<n$ we can consider the Riesz potential $I_\alpha$ defined by 
$$I_\alpha(f)(x)= \pi ^{-n/2}2^{-\alpha}\frac{\Gamma ((n-\alpha)/2)}{\Gamma (\alpha/2)}\int_{\mathbb{R}^n}\frac{f(y)}{|x-y|^{n-\alpha}}dy,$$
and we can define the operator $T_{\Omega, \alpha}$ by 
$$T_{\Omega,\alpha}(f)(x)=p.v.\int_{\mathbb{R}^n}\frac{\Omega(y/ |y|)}{|y|^{n+1-\alpha}}f(x-y)dy.$$
Thus, if $\Omega \in L^1(\mathbb{S}^{n-1})$, $\displaystyle{\int_{\mathbb{S}^{n-1}}\Omega\ d\sigma=0}$ and $\Omega \in L^{n,\infty}(\mathbb{S}^{n-1})$, in \cite[Theorem 1.1]{Hoang1} it was proven the following pointwise estimate
$$|T_{\Omega,\alpha}(f)(x)|\leq C\|\Omega\|_{L^{n,\infty}(\mathbb{S}^{n-1})}I_{\alpha}(|\nabla f|)(x), \qquad \mbox{for } f\in\mathcal{C}^\infty_0(\mathbb{R}^n).$$
More recently, in \cite{Hoang2}, a variant of these previous estimates with \emph{weighted} Riesz-type operators was considered and the study of weighted operators in pointwise inequalities as above will constitute the aim of this work. To this end, we recall now the definition of these weighted operators: if $f:\mathbb{R}^n\longrightarrow \mathbb{R}$ is a locally integrable function and if $\omega$ is a (generic) weight, then for a parameter $0<\alpha<n$, we define the weighted Riesz-type operator $\mathcal{F}_{\omega,\alpha}$ by the expression
$$\mathcal{F}_{\omega,\alpha}(f)(x)=\int_{\mathbb{R}^n}\frac{|x-y|^{\alpha}}{\omega(B(x,|x-y|))}f(y)\omega(y)dy.$$
Remark that in the absence of a weight (\emph{i.e.} replacing $\omega(B(x,|x-y|))$ by $|x-y|^n$ and $\omega(y)dy$ by $dy$), we recover (up to a constant) the classical Riesz potentials $I_\alpha$ defined above, indeed we have
$$\mathcal{F}_{\alpha}(f)(x)=\int_{\mathbb{R}^n}\frac{|x-y|^{\alpha}}{|x-y|^n}f(y)dy=\int_{\mathbb{R}^n}\frac{f(y)}{|x-y|^{n-\alpha}}dy\sim I_\alpha(f)(x).$$
With this ingredient, one of the main results of \cite{Hoang2} is the following: for the operator $T^*_\Omega$ defined in the expression (\ref{Def_OperatorTstar}) with $\Omega \in L^1(\mathbb{S}^{n-1})$, $\displaystyle{\int_{\mathbb{S}^{n-1}}\Omega \ d\sigma=0}$ and $\Omega\in L^{n,\infty}(\mathbb{S}^{n-1})$ and if $\omega$ is a Muckenhoupt weight that belongs to the $A_1$ class, then we have the estimate (see \cite[Theorem 2.3]{Hoang2})
\begin{equation}\label{EstimationAEtudier}
T^*_{\Omega}(f)(x)\leq C\|\Omega\|_{L^{n,\infty}(\mathbb{S}^{n-1})}[\omega]_{A_1}\mathcal{F}_{\omega,1}(|\nabla f|)(x) \qquad \mbox{for } f\in\mathcal{C}^\infty_0(\mathbb{R}^n).
\end{equation}
As pointed out in \cite[Remark 2.4]{Hoang2}, the previous estimate seems to be the first one to consider weighted Riesz potentials for this type of inequalities. Of course, interesting functional inequalities can be obtained from (\ref{EstimationAEtudier}) and in particular some weighted Sobolev like inequalities are obtained in \cite[Corollary 4.4.]{Hoang2}.\\

In this article we will generalize the estimate (\ref{EstimationAEtudier}) in two different directions. First we will consider a more general function $\Omega$ (in the sense studied in \cite{ChMarcociMarcoci2}) and second we will consider weights that belong to the wider Muckenhoupt class $A_\delta$ for some parameter $1<\delta<n$ (we have the class inclusion $A_1\subset A_\delta$, see \cite[Corollary 9.2.6]{Grafakos}). This framework will lead us to a more generic class of functional inequalities that will be studied below.\\ 

In order to state our main result, we need to introduce some notations. 
\begin{itemize}
\item[$\bullet$] {\bf Weighted maximal functions}. Consider $\omega:\mathbb{R}^n\longrightarrow \mathbb{R}^+$ be a weight. For a locally integrable function $f:\mathbb{R}^n\longrightarrow \mathbb{R}$, we define the weighted maximal function $\mathscr{M}_{\omega}$ of $f$ by the expression
\begin{equation}\label{Weightedmaximalfunctions}
\mathscr{M}_{\omega}(f)(x)=\displaystyle{\underset{B \ni x}{\sup } \;\frac{1}{\omega(B)}\int_{B }|f(y)|\omega(y)dy},
\end{equation}
where the supremum is taken over all open balls $B$ that contain the point $x$.

\item[$\bullet$] {\bf Weighted Morrey spaces}. For $1\leq p\leq q<+\infty$,  the weighted Morrey spaces $\mathcal{M}^{p,q}(\omega)$ are defined as the set of all measurable functions $f$ such that the condition
\begin{equation}\label{Def_Weighted_Morrey_space}
\|f\|_{\mathcal{M}^{p,q}(\omega)}=\underset{x\in \mathbb{R}^n, \; r>0}{\sup}\left(\frac{1}{\omega(B(x,r))^{1-\frac{p}{q}}}\int_{B(x,r)}|f(y)|^p \omega(y)dy\right)^{\frac{1}{p}}<+\infty
\end{equation}
is satisfied. Note that we have $L^q(\omega)=\mathcal{M}^{q,q}(\omega)$ and remark also that we have the space inclusions $\mathcal{M}^{p_1,q}(\omega)\subset \mathcal{M}^{p_0,q}(\omega)$ if $1\leq p_0\leq p_1\leq q$.
\item[$\bullet$] {\bf Lower Ahlfors condition for weights.} We will consider here Muckenhoupt weights $\omega\in A_\delta$ that satisfy the following \emph{${\bf d}$-lower Ahlfors condition} 
\begin{equation}\label{LowerAhlforsCondition}
Cr^{\bf d}\leq \omega(B(x,r)),\quad \mbox{for all } r>0, \  x\in \mathbb{R}^n,
\end{equation}
for some parameter ${\bf d}$ such that $1<{\bf d}<+\infty$

\end{itemize}
See Section \ref{Secc_FuncSpaces} below for more details on weighted maximal functions $\mathscr{M}_{\omega}$, weighted Morrey spaces $\mathcal{M}^{p,q}(\omega)$, Muckenhoupt weights, as well as some important properties between these objects.\\

Our first result is the following. 
\begin{Theoreme}[{\bf $A_\delta$-pointwise inequality}]\label{Theo_Pointwise_Weighted}
Consider the operator $T^*_\Omega$ defined in the expression (\ref{Def_OperatorTstar}) where the function $\Omega$ is such that $\Omega \in L^1(\mathbb{S}^{n-1})$, $\displaystyle{\int_{\mathbb{S}^{n-1}}\Omega \ d\sigma=0}$ and $\Omega\in L^\rho(\mathbb{S}^{n-1})$ with $1<\rho<n$.\\

\noindent Fix $\mathfrak{s}$ a real parameter such that $1<\overline{\rho}\leq \mathfrak{s} <n$ (where $\overline{\rho}=\frac{\rho n}{\rho n+\rho-n}$). Consider an $A_\delta$-weight $\omega$ (with $1\leq  \delta <+\infty$) that satisfies the ${\bf d}$-lower condition (\ref{LowerAhlforsCondition}) with ${\bf d}>\mathfrak{s}\delta$ and fix a parameter $\mathfrak{q}$ such that 
$$1<\mathfrak{s}\delta\leq \mathfrak{q}<{\bf d}.$$
\noindent If $f: \mathbb{R}^n \longrightarrow \mathbb{R}$ is a regular function such that $\nabla f$ belongs to the weighted Morrey space $\mathcal{M}^{\mathfrak{s}\delta, \mathfrak{q}}(\omega)$, then we have the following pointwise estimate
\begin{equation}\label{pointwise_estimate_p}
T^*_\Omega(f)(x)\leq C\|\Omega\|_{L^\rho(\mathbb{S}^{n-1})}[\omega]_{A_\delta}^\frac{1}{\mathfrak{s}\delta} \left(\mathscr{M}_\omega(|\nabla f|^{\mathfrak{s}\delta})(x)\right)^{\frac{1}{\mathfrak{s}\delta}(1-\frac{\mathfrak{q}}{{\bf d}})} \|\nabla f\|_{\mathcal{M}^{\mathfrak{s}\delta, \mathfrak{q}}(\omega)}^{\frac{\mathfrak{q}}{{\bf d}}}.
\end{equation}
\end{Theoreme}
\noindent As announced, this result is a double improvement of the estimate (\ref{EstimationAEtudier}) given in the article \cite{Hoang2}. Indeed, we can consider here a kernel $\Omega$ in the space $L^\rho(\mathbb{S}^{n-1})$ with $1<\rho<n$ and this constitutes a first improvement of the inequality (\ref{EstimationAEtudier}) since we have the space inclusions $L^{n,\infty}(\mathbb{S}^{n-1})\subset L^{\rho}(\mathbb{S}^{n-1})$ presented in the formula (\ref{SpaceInclusions0}) above. The second improvement is given by the fact that we can consider here a weight $\omega$ that belongs to the Muckenhoupt class $A_\delta$ with $1\leq \delta <+\infty$ (recall that we have the class inclusion $A_1\subset A_\delta$ for $1<\sigma<+\infty$). Note in particular that in \cite{Hoang2}, the only weights considered there were $\omega \in A_1$. 
Let us point out now that, since we are working over the entire space $\mathbb{R}^n$, the lower Ahlfors condition given in (\ref{LowerAhlforsCondition}) cannot be deduced from the properties of Muckenhoupt classes $A_\delta$, so this condition represents an additional constraint for the weights considered here (see Remark 4.2 in \cite{Hoang2} for more details on this particular point). Note finally  that the relationships $\mathfrak{s}<n$ and $\mathfrak{s}\delta\leq \mathfrak{q}<{\bf d}$ are essentially technical and we do not know if it is possible to get rid of these constraints or if it is possible to consider the cases $\mathfrak{s}\delta={\bf d}$ or even $\mathfrak{s}\delta>{\bf d}$. The study of these cases would probably require considering other approaches that are beyond the scope of this work.\\

Of course, the inequality (\ref{pointwise_estimate_p}) can be used as a starting point for several interesting estimates. In what follows we will always assume the following framework: we will consider $\Omega$ a function such that $\Omega \in L^1(\mathbb{S}^{n-1})$, $\displaystyle{\int_{\mathbb{S}^{n-1}}\Omega \ d\sigma=0}$ and such that $\Omega\in L^\rho(\mathbb{S}^{n-1})$ with $1<\rho<n$ and we will work with the operator $T^*_{\Omega}$ defined in (\ref{Def_OperatorTstar}).\\ 

We also recall that for a generic weight  $\omega$ and for $1\leq p<+\infty$, the weighted Lebesgue space $L^p(\omega)$ is given as the set of measurable functions such that 
$$\|f\|_{L^p(\omega)}=\left(\int_{\mathbb{R}^n}|f(x)|^p\omega(x)dx\right)^{\frac{1}{p}}<+\infty.$$
Under this general setting we have the following result:
\begin{Corollaire}[{\bf Weighted Sobolev estimates}]\label{Corollary2}
Consider a $A_\delta$-weight $\omega$, with $1\leq \delta<+\infty$, that satisfies the lower Ahlfors condition (\ref{LowerAhlforsCondition}) with index $1<{\bf d}<+\infty$. Fix $\mathfrak{s}$ a real parameter such that $1<\overline{\rho}=\frac{\rho n}{\rho n+\rho-n}\leq \mathfrak{s} <n$ and such that $1<\mathfrak{s}\delta<{\bf d}$.
\begin{itemize}
\item[1)] Consider a regular function $f:\mathbb{R}^n\longrightarrow \mathbb{R}$ such that $\|\nabla f\|_{L^{\mathfrak{q}}(\omega)}<+\infty$ with $1< q<+\infty$. If we have the relationship $\frac{r}{\mathfrak{s}\delta}(1-\frac{\mathfrak{q}}{\bf d})>1$ and if $r(1-\frac{\mathfrak{q}}{\bf d})=\mathfrak{q}$ for some index $1<r<+\infty$, then we have the following Sobolev-type inequality
$$\|T^*_\Omega(f)\|_{L^{r}(\omega)}\leq C_{\Omega, \omega}\|\nabla f\|_{L^{\mathfrak{q}}(\omega)}.$$
\item[2)] If $f:\mathbb{R}^n\longrightarrow \mathbb{R}$ is a regular function such that $\nabla f\in L^{\sigma}(\omega)$ with $1<\sigma<+\infty$ and $\nabla f\in \mathcal{M}^{\mathfrak{s}\delta, \mathfrak{q}}(\omega)$ with $\mathfrak{s}\delta\leq \mathfrak{q}<{\bf d}$. If $1<r<+\infty$ is a parameter such that $\frac{r}{\mathfrak{s}\delta}(1-\frac{\mathfrak{q}}{\bf d})>1$ and such that $r(1-\frac{\mathfrak{q}}{\bf d})=\sigma$, then we have the control
$$\|T^*_\Omega(f)\|_{L^r(\omega)}\leq C_{\Omega, \omega} \|\nabla f\|_{L^{\sigma}(\omega)}^{(1-\frac{\mathfrak{q}}{\bf d})}\|\nabla f\|_{\mathcal{M}^{\mathfrak{s}\delta, \mathfrak{q}}(\omega)}^{\frac{\mathfrak{q}}{\bf d}}.$$

\item[3)] Consider a parameter $1<\mathfrak{q}<+\infty$ such that $\mathfrak{s}\delta\leq \mathfrak{q}<{\bf d}$. Fix two indexes $1<\mathfrak{a}\leq \mathfrak{b}<+\infty$ such that $1< \mathfrak{a}({\frac{1}{\mathfrak{s}\delta}(1-\frac{\mathfrak{q}}{\bf d})})\leq  \mathfrak{b}({\frac{1}{\mathfrak{s}\delta}(1-\frac{\mathfrak{q}}{\bf d})})$. If $f:\mathbb{R}^n\longrightarrow \mathbb{R}$ is a regular function such that $\nabla f\in \mathcal{M}^{\mathfrak{a}(1-\frac{\mathfrak{q}}{\bf d}), \mathfrak{b}(1-\frac{\mathfrak{q}}{\bf d})}(\omega)$ and $\nabla f\in \mathcal{M}^{\mathfrak{s}\delta, \mathfrak{q}}(\omega)$, then we have the following functional inequality
$$\|T^*_\Omega(f)\|_{\mathcal{M}^{\mathfrak{a}, \mathfrak{b}}(\omega)}\leq  C_{\Omega, \omega} \|\nabla f\|_{\mathcal{M}^{\mathfrak{a}(1-\frac{\mathfrak{q}}{\bf d}), \mathfrak{b}(1-\frac{\mathfrak{q}}{\bf d})}(\omega)}^{(1-\frac{\mathfrak{q}}{\bf d})} \|\nabla f\|_{\mathcal{M}^{\mathfrak{s}\delta, \mathfrak{q}}(\omega)}^{\frac{\mathfrak{q}}{\bf d}}.$$
\end{itemize}
\end{Corollaire}
\noindent Weighted or unweighted functional inequalities in many different functional spaces were studied in \cite{Cianchi0}, \cite{Cianchi1}, \cite{Conde}, \cite{Derigoz}, \cite{GCRdF}, \cite{Hoang}, \cite{Hoang1}, \cite{Hoang2}, \cite{Kinnunen}, \cite{Komori}, \cite{Krbec}, \cite{Li}, \cite{RaSam}. The Sobolev-like inequalities above are direct consequences of the previous pointwise estimates (\ref{pointwise_estimate_p}) and, to the best of our knowledge, these are new functional inequalities. Let us note in particular that the last inequality of the previous corollary can be viewed as an interpolation result in the following sense: if we define the space $\mathbb{M}^{\mathfrak{a}, \mathfrak{b}}_{T^*_\Omega}(\omega)$ as the set of measurable functions $f:\mathbb{R}^n\longrightarrow \mathbb{R}$ such that $\|f\|_{\mathbb{M}^{\mathfrak{a}, \mathfrak{b}}_{T^*_\Omega}(\omega)}=\|T^*_\Omega(f)\|_{\mathcal{M}^{\mathfrak{a}, \mathfrak{b}}(\omega)}<+\infty$, where we have $1<\mathfrak{a}\leq \mathfrak{b}<+\infty$, and if we define the Sobolev like space $\dot{\mathcal{W}}^1_{\mathcal{M}^{\mathfrak{a}, \mathfrak{b}}}(\omega)$ by the condition $\|f\|_{\dot{\mathcal{W}}^1_{\mathcal{M}^{\mathfrak{a}, \mathfrak{b}}}(\omega)}=\|\nabla f\|_{\mathcal{M}^{\mathfrak{a}, \mathfrak{b}}(\omega)}<+\infty$, then we can write
$$\left[\dot{\mathcal{W}}^1_{\mathcal{M}^{\mathfrak{a}_0, \mathfrak{b}_0}(\omega)}, \dot{\mathcal{W}}^1_{\mathcal{M}^{\mathfrak{a}_1, \mathfrak{b}_1}(\omega)}\right]_{\theta}=\mathbb{M}^{\mathfrak{a}, \mathfrak{b}}_{T^*_\Omega}(\omega),$$
for some parameter $0<\theta<1$ and where of course the indexes $\mathfrak{a}_0, \mathfrak{b}_0, \mathfrak{a}_1,\mathfrak{a}_1, \mathfrak{a}, \mathfrak{b}$ satisfy the conditions stated in the last 
point of the corollary above. In this weighted context the interpolation relationship stated above seems to be new, however a systematic study of these interpolation results falls beyond the scope of this article. See the books \cite{BS} and \cite{Bergh} for more details on interpolation theory.\\

The plan of the article is the following. In Section \ref{Secc_FuncSpaces} we will recall the definitions and the main properties of the functional spaces used here in connection with Muckenhoupt weights. In Section \ref{Secc_ProofTheo_Pointwise_Weighted} we will prove Theorem \ref{Theo_Pointwise_Weighted}. Next, in Section \ref{Secc_ProofCorollaries} we will prove the Corollary \ref{Corollary2}. In this paper the constants may vary at each occurrence.
\mysection{Some functional spaces, classical inequalities and Muckenhoupt weights properties}\label{Secc_FuncSpaces}
In this section we recall the definitions and some well known properties of the functional spaces, Muckenhoupt weights and maximal function operators that will be used here. 
\subsubsection*{Unweighted Functional spaces}
\begin{itemize}
\item For $1\leq p<+\infty$ and for $A=\mathbb{R}^n$ or $A \subset \mathbb{R}^n$, the usual Lebesgue space $L^p(A)$ are defined by the classical condition
$\|f\|_{L^p}=\displaystyle{\left(\int_{A}|f(x)|^pdx\right)^{\frac{1}{p}}}<+\infty$. Recall in particular that if $A$ is a bounded subset then we have the space inclusions $L^{p_1}(A)\subset L^{p_0}(A)\subset L^1(A)$ for $1\leq p_0\leq p_1$. Of course these inclusions are still valid if we consider $A=\mathbb{S}^{n-1}$.

\item For $1\leq p<+\infty$, Lorentz spaces $L^{p,\infty}(A)$ with $A=\mathbb{R}^n$ or $A=\mathbb{S}^{n-1}$ are defined by the condition $\|f\|_ {L^{p,\infty}}=\underset{\lambda>0}{\sup}\{\lambda \times |\{x\in A:|f(x)|> \lambda\}|^{1/p}\}<+\infty$. Recall now that by the real interpolation theory (see \cite[Theorem 5.2.1]{Bergh}) we have for some parameter $0<\theta<1$ the identity
$$(L^p(A), L^\infty(A))_{\theta,\infty}=L^{\frac{p}{1-\theta},\infty}(A).$$
Recall that we always have $L^{\frac{p}{1-\theta},\infty}(A)\subset L^p(A)+L^\infty(A)$, but if the set $A\subset \mathbb{R}^n$ is bounded, we also have the space inclusions
$$L^{\frac{p}{1-\theta},\infty}(A)\subset L^p(A)+L^\infty(A)\subset L^p(A)+L^p(A)=L^p(A).$$ 
Thus, in the particular case of $A=\mathbb{S}^{n-1}$, since $\sigma(\mathbb{S}^{n-1})<+\infty$, we deduce that the Lorentz spaces $L^{q,\infty}(\mathbb{S}^{n-1})$ are embedded in the Lebesgue spaces $L^\rho(\mathbb{S}^{n-1})$ as long as $q>\rho$. In particular, we have 
\begin{equation}\label{InclusionLorentzLpUloc}
L^{n,\infty}(\mathbb{S}^{n-1})\subset L^\rho(\mathbb{S}^{n-1}), \qquad \mbox{if }\quad 1\leq \rho<n.
\end{equation}
\end{itemize}
\subsubsection*{One classical inequality}
We need to recall an important inequality. Indeed, for a function $f\in \mathcal{C}^\infty_0(\mathbb{R}^n)$ and for all ball $B(x,r)$ such that $B(x,r)\subset supp(f)$ we have the following \emph{Poincaré-Sobolev inequality}:
\begin{equation}\label{PoincareSobolev_inequality}
\left(\frac{1}{|B(x,r)|}\int_{B(x,r)}|f(y)-f_{B_r}|^{q}dy\right)^{\frac{1}{q}}\leq Cr\left(\frac{1}{|B(x,r)|}\int_{B(x,r)}|\nabla f(y)|^{\mathfrak{s}}dy\right)^{\frac{1}{\mathfrak{s}}}. 
\end{equation}
for $1\leq \mathfrak{s}<n$ and $1\leq q\leq \frac{n\mathfrak{s}}{n-\mathfrak{s}}$. 
See a proof of this inequality in \cite[Theorem 3.14]{Kinnunen}. 
\subsubsection*{Properties of the Muckenhoupt weights}
We start recalling that a generic weight $\omega:\mathbb{R}^n\longrightarrow ]0,+\infty[$ is a locally integrable function. For a given weight $\omega$ and a measurable set $A$, we denote by $\omega(A)=\displaystyle{\int_A} \omega(x) \ dx$. We say that a weight $\omega$ belongs to the Muckenhoupt class $A_{\delta}$ for some $1<\delta<+\infty$ if we have
$$[\omega]_{A_{\delta}}=\underset{B}{\sup}\left(\frac{1}{|B|}\int_{B}\omega(x)dx\right)\left(\frac{1}{|B|}\int_{B}\omega(x)^{-\frac{1}{\delta-1}}dx\right)^{\delta-1}<+\infty,$$
where the supremum is taken over all balls $B$ in $\mathbb{R}^n$ and we will say that $\omega\in A_1$ if 
$$[\omega]_{A_{1}}=\underset{B}{\sup}\left(\frac{1}{|B|}\int_{B}\omega(x)dx\right)\|\omega^{-1}\|_{L^\infty(B)}<+\infty,$$
where the supremum is taken over all balls $B$ in $\mathbb{R}^n$.\\

From this we easily deduce that, if $A\subset B$ are bounded measurable sets of $\mathbb{R}^n$, we then have $\omega(A)\leq \omega(B)$. Muckenhoupt weights have several properties, in particular, if  $\omega\in A_{\delta}$ with $1\leq \delta<+\infty$, then we have the doubling property for $\omega$, i.e. there exists a constant $C=C(\omega)>0$ such that
\begin{equation}\label{DoublingProperty}
\omega(B(x,2r))\leq C \omega(B(x,r)), \quad \mbox{for all } x\in \mathbb{R}^n,  r>0.
\end{equation}
The previous property can also be expressed in the following manner
\begin{equation}\label{DoublingProperty1}
\omega(B(x,\lambda r))\leq \lambda^{n \delta}[\omega]_{A_{\delta}}\omega(B(x,r)),
\end{equation}
for all $\lambda>1$, $x\in \mathbb{R}^n$ and $r>0$.
See point (9) of Proposition 9.1.5 of \cite{Grafakos} for more details about this doubling property.
Next by using H\" older's inequality we have the following estimate:

\begin{Lemme}\label{Lemme_Aweights}
If $\omega$ is a $A_{\delta}$ weight with $1\leq \delta<+\infty$, then by using H\" older's inequality we have the following estimate
$$\frac{1}{|B|}\int_{B}|f(x)|dx\leq [\omega]_{A_{\delta}}^\frac{1}{{\delta}}\left(\frac{1}{\omega(B)}\int_B|f(x)|^{\delta}\omega(x)dx\right)^\frac{1}{{\delta}},$$
for any ball $B\subset \mathbb{R}^n$.
\end{Lemme}
\noindent See the book \cite{Grafakos}, p. 285, for a proof of this fact as well as for more properties of the Muckenhoupt classes $A_{\delta}$.
\subsubsection*{Weighted functional spaces}
We need now to consider weighted functional spaces and to recall some of their properties. 
\begin{itemize}
\item {\bf Weighted Lebesgue spaces}. For $1\leq p<+\infty$ and for $A=\mathbb{R}^n$ or $A \subset \mathbb{R}^n$, the weighted Lebesgue spaces $L^{p}(\omega)$ are defined by the condition
$$\|f\|_{L^{p}(\omega)}=\left(\int_{A}|f(x)|^{p}\omega(x)dx\right)^{\frac{1}{{p}}}<+\infty.$$
Note that for some parameter $1<\alpha<+\infty$ we have the property 
\begin{equation}\label{ProprietePuissanceLebesgue}
\||f|^\alpha\|_{L^p(\omega)}=\|f\|^\alpha_{L^{\alpha p}(\omega)}.
\end{equation}
Remark also that in this weighted setting we have the following H\"older inequality for $\frac{1}{q}+\frac{1}{q'}=1$, for some measurable set $A\subseteq \mathbb{R}^n$ and for $\omega\in A_\delta$ with $1\leq \delta<+\infty$:
$$\int_{A}|f(x)|\omega(x)dx\leq \|f\|_{L^q_\omega(A)}\|f\|_{L^{q'}_\omega(A)}.$$
We recall now that the weighted maximal function $\mathscr{M}_\omega$ given in (\ref{Weightedmaximalfunctions}) is bounded in the weighted Lebesgue spaces and we have 
$$\|\mathscr{M}_\omega(f)\|_{L^p(\omega)}\leq C\|f\|_{L^p(\omega)},$$
where $1<p<+\infty$ and $\omega\in A_\delta$ for $1\leq \delta <+\infty$. See a proof of this fact in \cite[Theorem 2.6, p. 146]{GCRdF}.
\item {\bf Weighted Morrey spaces}. For $1<p\leq q<+\infty$ we recall the definition of the weighted Morrey spaces, which are given by the condition $\|f\|_{\mathcal{M}^{p,q}(\omega)}<+\infty$ where
$$\|f\|_{\mathcal{M}^{p,q}(\omega)}=\underset{x\in \mathbb{R}^n, \; r>0}{\sup}\left(\frac{1}{\omega(B(x,r))^{1-\frac{p}{q}}}\int_{B(x,r)}|f(y)|^p \omega(y)dy\right)^{\frac{1}{p}}<+\infty.$$
A straightforward but useful property is the following:
\begin{equation}\label{ProprietePuissanceMorrey}
\||f|^\alpha\|_{\mathcal{M}^{p,q}(\omega)}=\|f\|_{\mathcal{M}^{\alpha p, \alpha q}(\omega)}^\alpha.
\end{equation}
\end{itemize}
Following \cite[Theorem 3.1]{Komori}, the maximal function $\mathscr{M}_\omega$ is also bounded in these weighted Morrey spaces and we have 
$$\|\mathscr{M}_\omega(f)\|_{\mathcal{M}^{p,q}(\omega)}\leq C\|f\|_{\mathcal{M}^{p,q}(\omega)},$$
where $\omega \in A_{\delta}$ for $1\leq \delta<+\infty$ and where $1<p\leq q<+\infty$.
\mysection{Proof of the Theorem \ref{Theo_Pointwise_Weighted}}\label{Secc_ProofTheo_Pointwise_Weighted}
Let us start by recalling the definition of the operator $T^*_{\Omega}$ given in the expression  (\ref{Def_OperatorTstar}) above: we thus have
$$T^*_{\Omega}(f)(x)=\underset{t>0}{\sup}\left|\int_{\{|y|>t\}}\frac{\Omega(y/|y|)}{|y|^n}f(x-y)dy\right|,$$
and for some $t>0$ we consider now the operator
$$T^t_{\Omega}(f)(x)=\int_{\{|y|>t\}}\frac{\Omega(y/|y|)}{|y|^n}f(x-y)dy.$$
By definition, we easily note that $T^*_{\Omega}(f)(x)=\underset{t>0}{\sup}|T^t_{\Omega}(f)(x)|$ and that $|T_\Omega(f)(x)|\leq T^*_{\Omega}(f)(x)$.\\

Now, for a function $f\in \mathcal{C}^\infty_0(\mathbb{R}^n)$ and for some $k_0\in \mathbb{Z}$ so that $2^{k_0-2}<t \leq 2^{k_0-1}$, we write
$$T^t_{\Omega}(f)(x)=\int_{\{t<|y|\leq 2^{k_0-1}\}}\frac{\Omega(y/|y|)}{|y|^n}f(x-y)dy+\sum_{k\geq k_0}\int_{\{2^{k-1}<|y|\leq 2^{k}\}}\frac{\Omega(y/|y|)}{|y|^n}f(x-y)dy.$$
Using the fact that the function $\Omega$ is of null integral, we can introduce some constants in the previous expression to obtain
$$T^t_{\Omega}(f)(x)=\int_{\{t<|y|\leq 2^{k_0-1}\}}\frac{\Omega(y/|y|)}{|y|^n}(f(x-y)-c_{k_0})dy+\sum_{k\geq k_0}\int_{\{2^{k-1}<|y|\leq 2^{k}\}}\frac{\Omega(y/|y|)}{|y|^n}(f(x-y)-c_k)dy,$$
from which we deduce
\begin{eqnarray*}
|T^t_{\Omega}(f)(x)|&\leq &\sum_{k\in \mathbb{Z}}\int_{\{2^{k-1}<|y|\leq 2^{k}\}}|\frac{\Omega(y/|y|)}{|y|^n}| |f(x-y)-c_k|dy\\
&\leq &C\sum_{k\in \mathbb{Z}}\frac{1}{2^{kn}}\int_{\{|y|\leq 2^{k}\}}|\Omega(y/|y|)| |f(x-y)-c_k|dy,
\end{eqnarray*}
Now, by the H\"older inequality with $\frac{1}{\rho}+\frac{1}{\rho'}=1$ and $1<\rho<n$, we write
$$|T^t_{\Omega}(f)(x)|\leq C\sum_{k\in \mathbb{Z}}\frac{1}{2^{kn}}\left(\int_{\{|y|\leq 2^{k}\}}|\Omega(y/|y|)|^\rho dy\right)^{\frac{1}{\rho}}\left(\int_{\{|y|\leq 2^{k}\}}|f(x-y)-c_k|^{\rho'}dy\right)^{\frac{1}{\rho'}}.$$
Introducing the variable $z=2^{-k} y$, by a change of variables in the first integral we obtain
$$|T^t_{\Omega}(f)(x)|\leq C\sum_{k\in \mathbb{Z}}\frac{1}{2^{kn(1-\frac{1}{\rho})}}\left(\int_{\{|z|\leq 1\}}|\Omega(z/|z|)|^\rho dz\right)^{\frac{1}{\rho}}\left(\int_{\{|y|\leq 2^{k}\}}|f(x-y)-c_k|^{\rho'}dy\right)^{\frac{1}{\rho'}},$$
and rewriting this formula we have
\begin{eqnarray*}
|T^t_{\Omega}(f)(x)|&\leq &C\sum_{k\in \mathbb{Z}}\left(\int_{\{|z|\leq 1\}}|\Omega(z/|z|)|^\rho dz\right)^{\frac{1}{\rho}}\frac{1}{2^{\frac{kn}{\rho'}}}\left(\int_{\{|y|\leq 2^{k}\}}|f(x-y)-c_k|^{\rho'}dy\right)^{\frac{1}{\rho'}}\\
&\leq &C\sum_{k\in \mathbb{Z}}\left(\int_{\{|z|\leq 1\}}|\Omega(z/|z|)|^\rho dz\right)^{\frac{1}{\rho}}\left(\frac{1}{2^{kn}}\int_{\{|y|\leq 2^{k}\}}|f(x-y)-c_k|^{\rho'}dy\right)^{\frac{1}{\rho'}}.
\end{eqnarray*}
For the second integral above, we consider now the ball $B(x,2^k)$ and we fix the constant $c_k=f_{B_k}=\frac{1}{|B(x,2^k)|}\displaystyle{\int_{B(x,2^k)}f(y)dy}$, so we can write (since $v_n 2^{kn}=|B(x,2^k)|$, where $v_n=|B(0,1)|$ is the volume of the $n$-dimensional unit ball):
$$|T^t_{\Omega}(f)(x)|\leq C\sum_{k\in \mathbb{Z}}\left(\int_{\{|z|\leq 1\}}|\Omega(z/|z|)|^\rho dz\right)^{\frac{1}{\rho}}\left(\frac{1}{|B(x,2^k)|}\int_{B(x,2^k)}|f(y)-f_{B_k}|^{\rho'}dy\right)^{\frac{1}{\rho'}}.$$
We study now more in detail the first integral above, we thus have
\begin{eqnarray*}
|T^t_{\Omega}(f)(x)|&\leq &C\left(\int_{0}^1\int_{\mathbb{S}^{n-1}}|\Omega(\xi/|\xi|)|^\rho d\sigma(\xi)r^{n-1}dr\right)^{\frac{1}{\rho}}\sum_{k\in \mathbb{Z}}\left(\frac{1}{|B(x,2^k)|}\int_{B(x,2^k)}|f(y)-f_{B_k}|^{\rho'}dy\right)^{\frac{1}{\rho'}}\\
&\leq &C\left(\int_{\mathbb{S}^{n-1}}|\Omega(\xi/|\xi|)|^\rho d\sigma(\xi)\right)^{\frac{1}{\rho}}\sum_{k\in \mathbb{Z}}\left(\frac{1}{|B(x,2^k)|}\int_{B(x,2^k)}|f(y)-f_{B_k}|^{\rho'}dy\right)^{\frac{1}{\rho'}},
\end{eqnarray*}
so we obtain
$$|T^t_{\Omega}(f)(x)|\leq C \|\Omega\|_{L^\rho(\mathbb{S}^{n-1})}\sum_{k\in \mathbb{Z}}\left(\frac{1}{|B(x,2^k)|}\int_{B(x,2^k)}|f(y)-f_{B_k}|^{\rho'}dy\right)^{\frac{1}{\rho'}}.$$
Note that since $1< \rho<n$, by (\ref{InclusionLorentzLpUloc}) we have $L^{n,\infty}(\mathbb{S}^{n-1})\subset L^\rho(\mathbb{S}^{n-1})$ and thus the norm $ \|\Omega\|_{L^\rho(\mathbb{S}^{n-1})}$ induces a refinement with respect to the norm $ \|\Omega\|_{L^{n,\infty}}$ used in \cite{Hoang} and in \cite{Hoang2}. Note also that if $\rho<n$ then we have $\frac{n}{n-1}<\rho'$ as we have $\frac{1}{\rho}+\frac{1}{\rho'}=1$.\\

\noindent Now we apply the Poincaré-Sobolev inequality given in (\ref{PoincareSobolev_inequality}) to obtain
\begin{eqnarray}
|T^t_{\Omega}(f)(x)|&\leq &C \|\Omega\|_{L^\rho(\mathbb{S}^{n-1})}\sum_{k\in \mathbb{Z}}\left(\frac{1}{|B(x,2^k)|}\int_{B(x,2^k)}|f(y)-f_{B_k}|^{\rho'}dy\right)^{\frac{1}{\rho'}}\notag\\
&\leq &C \|\Omega\|_{L^\rho(\mathbb{S}^{n-1})}\underbrace{\sum_{k\in \mathbb{Z}}2^k\left(\frac{1}{|B(x,2^k)|}\int_{B(x,2^k)}|\nabla f(y)|^{\mathfrak{s}}dy\right)^{\frac{1}{\mathfrak{s}}}}_{\mathcal{S}},\label{SumPointWiseInequality}
\end{eqnarray}
where $\frac{n}{n-1}<\rho'$ (since $1<\rho<n$ and $\frac{1}{\rho}+\frac{1}{\rho'}=1$) and where $\rho'\leq \frac{n\mathfrak{s}}{n-\mathfrak{s}}$. Note that we thus have $\frac{n}{n-1}<\rho'\leq \frac{n\mathfrak{s}}{n-\mathfrak{s}}$ which leads us to the condition $1<\frac{\rho n}{\rho n+\rho-n}=\overline{\rho}\leq \mathfrak{s}<n$. \\

\noindent We study now the sum $\mathcal{S}$ in the previous formula and we have
$$\mathcal{S}=\sum_{k\in \mathbb{Z}}2^k\left(\frac{1}{|B(x,2^k)|}\int_{B(x,2^k)}|\nabla f(y)|^{\mathfrak{s}}dy\right)^{\frac{1}{\mathfrak{s}}},$$
which we rewrite as follows
\begin{equation}\label{FormuleAB}
\begin{split}
\mathcal{S}&\leq \underbrace{\sum_{k\in \mathbb{Z}}2^k\left(\frac{1}{|B(x,2^k)|}\int_{\{2^{k-1}<|x-y|\leq 2^{k}\}}|\nabla f(y)|^{\mathfrak{s}}dy\right)^{\frac{1}{\mathfrak{s}}}}_{(A)}\\
&+\underbrace{\sum_{k\in \mathbb{Z}}2^k\left(\frac{1}{|B(x,2^k)|}\int_{\{|x-y|\leq 2^{k-1}\}}|\nabla f(y)|^{\mathfrak{s}}dy\right)^{\frac{1}{\mathfrak{s}}}}_{(B)}.
\end{split}
\end{equation}
We will estimate each of the previous terms separately. 
\begin{itemize}
\item For the term $(A)$ above we write
$$(A)=\sum_{k\in \mathbb{Z}}2^k\left(\frac{1}{|B(x,2^k)|}\int_{\{2^{k-1}<|x-y|\leq 2^{k}\}}|\nabla f(y)|^{\mathfrak{s}}dy\right)^{\frac{1}{\mathfrak{s}}},$$
and introducing a technical parameter $0<\mathcal{K}<+\infty$ that will be fixed later we obtain
\begin{equation}\label{SommeA}
\begin{split}
(A)&\leq \underbrace{\sum_{k\leq \lfloor \log_2(\mathcal{K})\rfloor}2^k\left(\frac{1}{|B(x,2^k)|}\int_{B(x,2^k)}|\nabla f(y)|^{\mathfrak{s}}dy\right)^{\frac{1}{\mathfrak{s}}}}_{(A_1)}\\
&+ \underbrace{\sum_{k> \lfloor \log_2(\mathcal{K})\rfloor}2^k\left(\frac{1}{|B(x,2^k)|}\int_{B(x,2^k)}|\nabla f(y)|^{\mathfrak{s}}dy\right)^{\frac{1}{\mathfrak{s}}}}_{(A_2)}.
\end{split}
\end{equation}
For the term $(A_1)$ above, by applying the Lemma \ref{Lemme_Aweights} we obtain
\begin{eqnarray*}
(A_1)&=&\sum_{k\leq \lfloor \log_2(\mathcal{K})\rfloor}2^k\left(\frac{1}{|B(x,2^k)|}\int_{B(x,2^k)}|\nabla f(y)|^{\mathfrak{s}}dy\right)^{\frac{1}{\mathfrak{s}}}\\
&\leq &\sum_{k\leq \lfloor \log_2(\mathcal{K})\rfloor}2^k[\omega]_{A_\delta}^\frac{1}{\mathfrak{s}\delta}\left(\frac{1}{\omega(B(x,2^k))}\int_{B(x,2^k)}|\nabla f(y)|^{\mathfrak{s}\delta}\omega(y)dy\right)^{\frac{1}{\mathfrak{s}\delta}}.
\end{eqnarray*}
Now, using the definition of the weighted maximal functions given in the expression (\ref{Weightedmaximalfunctions}) above, we have
\begin{eqnarray}
(A_1)&\leq &[\omega]_{A_\delta}^\frac{1}{\mathfrak{s}\delta} \left(\mathscr{M}_\omega(|\nabla f|^{\mathfrak{s}\delta})(x)\right)^{\frac{1}{\mathfrak{s}\delta}}\sum_{k\leq \lfloor \log_2(\mathcal{K})\rfloor}2^k\notag\\
&\leq &C [\omega]_{A_\delta}^\frac{1}{\mathfrak{s}\delta} \left(\mathscr{M}_\omega(|\nabla f|^{\mathfrak{s}\delta})(x)\right)^{\frac{1}{\mathfrak{s}\delta}} \mathcal{K}.\label{FormuleA1}
\end{eqnarray}
For the term $(A_2)$ of the expression (\ref{SommeA}), we write
$$(A_2)=\sum_{k> \lfloor \log_2(\mathcal{K})\rfloor}2^k\left(\frac{1}{|B(x,2^k)|}\int_{B(x,2^k)}|\nabla f(y)|^{\mathfrak{s}}dy\right)^{\frac{1}{\mathfrak{s}}}.$$
Applying the Lemma \ref{Lemme_Aweights} again, we have the inequality
$$(A_2)\leq \sum_{k> \lfloor \log_2(\mathcal{K})\rfloor}2^k[\omega]_{A_\delta}^{\frac{1}{\mathfrak{s}\delta}}\left(\frac{1}{\omega(B(x,2^k))}\int_{B(x,2^k)}|\nabla f(y)|^{\mathfrak{s}\delta}\omega(y)dy\right)^{\frac{1}{\mathfrak{s}\delta}},$$
which we rewrite as follows
$$(A_2)\leq [\omega]_{A_\delta}^{\frac{1}{\mathfrak{s}\delta}}\sum_{k> \lfloor \log_2(\mathcal{K})\rfloor}2^k\left(\frac{\omega(B(x,2^k))^{-\frac{\mathfrak{s}\delta}{\mathfrak{q}}}}{\omega(B(x,2^k))^{1-\frac{\mathfrak{s}\delta}{\mathfrak{q}}}}\int_{B(x,2^k)}|\nabla f(y)|^{\mathfrak{s}\delta}\omega(y)dy\right)^{\frac{1}{\mathfrak{s}\delta}}.$$
At this point we use the Ahlfors condition $Cr^{\bf d}\leq \omega(B(x,r))$ assumed by hypothesis (\emph{i.e.} $\omega(B(x,2^k))^{-\frac{\mathfrak{s}\delta}{\mathfrak{q}}}\leq C(2^k)^{-{\bf d}\frac{\mathfrak{s}\delta}{\mathfrak{q}}}$) and we obtain the inequality
\begin{eqnarray*}
(A_2)&\leq &C [\omega]_{A_\delta}^{\frac{1}{\mathfrak{s}\delta}}\sum_{k> \lfloor \log_2(\mathcal{K})\rfloor}2^k\left(\frac{(2^k)^{-{\bf d}\frac{\mathfrak{s}\delta}{\mathfrak{q}}}}{\omega(B(x,2^k))^{1-\frac{\mathfrak{s}\delta}{\mathfrak{q}}}}\int_{B(x,2^k)}|\nabla f(y)|^{\mathfrak{s}\delta}\omega(y)dy\right)^{\frac{1}{\mathfrak{s}\delta}}\\
&\leq & C[\omega]_{A_\delta}^{\frac{1}{\mathfrak{s}\delta}}\sum_{k> \lfloor \log_2(\mathcal{K})\rfloor}2^{k(1-\frac{\bf d}{\mathfrak{q}})}\left(\frac{1}{\omega(B(x,2^k))^{1-\frac{\mathfrak{s}\delta}{\mathfrak{q}}}}\int_{B(x,2^k)}|\nabla f(y)|^{\mathfrak{s}\delta}\omega(y)dy\right)^{\frac{1}{\mathfrak{s}\delta}}.
\end{eqnarray*}
Now, using the definition of the weighted Morrey spaces $\mathcal{M}^{p,q}(\omega)$ given in (\ref{Def_Weighted_Morrey_space}) we can write
$$(A_2)\leq  C[\omega]_{A_\delta}^{\frac{1}{\mathfrak{s}\delta}}\sum_{k> \lfloor \log_2(\mathcal{K})\rfloor}2^{k(1-\frac{\bf d}{\mathfrak{q}})}\|\nabla f\|_{\mathcal{M}^{\mathfrak{s}\delta, \mathfrak{q}}(\omega)}.$$
Since by hypothesis we have $\mathfrak{s}\delta\leq \mathfrak{q}<{\bf d}$, the previous sum converges and we obtain
\begin{equation}\label{FormuleA2}
(A_2)\leq  C[\omega]_{A_\delta}^{\frac{1}{\mathfrak{s}\delta}}\|\nabla f\|_{\mathcal{M}^{\mathfrak{s}\delta, \mathfrak{q}}(\omega)} \mathcal{K}^{1-\frac{\bf d}{\mathfrak{q}}}.
\end{equation}
With the estimates (\ref{FormuleA1}) for the term $(A_1)$ and with the previous control (\ref{FormuleA2}) for $(A_2)$, getting back to the inequality (\ref{SommeA}), we obtain
\begin{eqnarray*}
(A)&\leq &(A_1)+(A_2)\\
&\leq & C [\omega]_{A_\delta}^\frac{1}{\mathfrak{s}\delta} \left(\mathscr{M}_\omega(|\nabla f|^{\mathfrak{s}\delta})(x)\right)^{\frac{1}{\mathfrak{s}\delta}} \mathcal{K}+C[\omega]_{A_\delta}^{\frac{1}{s\delta}}\|\nabla f\|_{\mathcal{M}^{\mathfrak{s}\delta, \mathfrak{q}}(\omega)} \mathcal{K}^{1-\frac{d}{\mathfrak{q}}}.
\end{eqnarray*}
Now we set $\mathcal{K}=\left(\frac{\|\nabla f\|_{\mathcal{M}^{\mathfrak{s}\delta, \mathfrak{q}}(\omega)}}{\left(\mathscr{M}_\omega(|\nabla f|^{\mathfrak{s}\delta})(x)\right)^{\frac{1}{\mathfrak{s}\delta}} }\right)^{\frac{\mathfrak{q}}{{\bf d}}}$ and we obtain
\begin{equation}\label{Estimation_Terme_A}
(A)\leq  C [\omega]_{A_\delta}^\frac{1}{\mathfrak{s}\delta} \left(\mathscr{M}_\omega(|\nabla f|^{\mathfrak{s}\delta})(x)\right)^{\frac{1}{\mathfrak{s}\delta}(1-\frac{\mathfrak{q}}{\bf d})} \|\nabla f\|_{\mathcal{M}^{\mathfrak{s}\delta, \mathfrak{q}}(\omega)}^{\frac{\mathfrak{q}}{\bf d}}.
\end{equation}
\item For the term (B) in (\ref{FormuleAB}) we have:
\begin{eqnarray*}
(B)&=&\sum_{k\in \mathbb{Z}}2^k\left(\frac{1}{|B(x,2^k)|}\int_{\{|x-y|\leq 2^{k-1}\}}|\nabla f(y)|^{\mathfrak{s}}dy\right)^{\frac{1}{\mathfrak{s}}}\\
&=&2^{1-\frac{n}{\mathfrak{s}}}\sum_{k\in \mathbb{Z}}2^{k-1}\left(\frac{1}{|B(x,2^{k-1})|}\int_{\{|x-y|\leq 2^{k-1}\}}|\nabla f(y)|^{\mathfrak{s}}dy\right)^{\frac{1}{\mathfrak{s}}},
\end{eqnarray*}
and we can thus write
\begin{equation}\label{Estimation_Terme_B}
(B)\leq  2^{1-\frac{n}{\mathfrak{s}}}\sum_{k\in \mathbb{Z}}2^{k}\left(\frac{1}{|B(x,2^k)|}\int_{\{|x-y|\leq 2^{k}\}}|\nabla f(y)|^{\mathfrak{s}}dy\right)^{\frac{1}{\mathfrak{s}}}= 2^{1-\frac{n}{\mathfrak{s}}}\mathcal{S}.
\end{equation}
\end{itemize}
With the estimates (\ref{Estimation_Terme_A}) and (\ref{Estimation_Terme_B}) at our disposal, coming back to the expression (\ref{FormuleAB}) we obtain
$$\mathcal{S}\leq  C [\omega]_{A_\delta}^\frac{1}{\mathfrak{s}\delta} \left(\mathscr{M}_\omega(|\nabla f|^{\mathfrak{s}\delta})(x)\right)^{\frac{1}{\mathfrak{s}\delta}(1-\frac{\mathfrak{q}}{\bf d})}  \|\nabla f\|_{\mathcal{M}^{\mathfrak{s}\delta, \mathfrak{q}}(\omega)}^{\frac{\mathfrak{q}}{\bf d}}+2^{1-\frac{n}{\mathfrak{s}}}\mathcal{S},$$
which is 
$$\mathcal{S}(1-2^{1-\frac{n}{\mathfrak{s}}})\leq  C [\omega]_{A_\delta}^\frac{1}{\mathfrak{s}\delta} \left(\mathscr{M}_\omega(|\nabla f|^{\mathfrak{s}\delta})(x)\right)^{\frac{1}{\mathfrak{s}\delta}(1-\frac{\mathfrak{q}}{\bf d})} \|\nabla f\|_{\mathcal{M}^{\mathfrak{s}\delta, \mathfrak{q}}(\omega)}^{\frac{\mathfrak{q}}{\bf d}},$$
since $\mathfrak{s}<n$, we have $2^{1-\frac{n}{\mathfrak{s}}}<1$ and we can write
$$\mathcal{S}\leq  \frac{C}{(1-2^{1-\frac{n}{\mathfrak{s}}})} [\omega]_{A_\delta}^\frac{1}{\mathfrak{s}\delta} \left(\mathscr{M}_\omega(|\nabla f|^{\mathfrak{s}\delta})(x)\right)^{\frac{1}{\mathfrak{s}\delta}(1-\frac{\mathfrak{q}}{\bf d})}  \|\nabla f\|_{\mathcal{M}^{\mathfrak{s}\delta, \mathfrak{q}}(\omega)}^{\frac{\mathfrak{q}}{\bf d}}.$$
With this estimate for the quantity $\mathcal{S}$, we can come back to the inequality (\ref{SumPointWiseInequality}) to obtain the control 
\begin{eqnarray*}
|T^t_{\Omega}(f)(x)|&\leq &C \|\Omega\|_{L^\rho(\mathbb{S}^{n-1})}\underbrace{\sum_{k\in \mathbb{Z}}2^k\left(\frac{1}{|B(x,2^k)|}\int_{B(x,2^k)}|\nabla f(y)|^{\mathfrak{s}}dy\right)^{\frac{1}{\mathfrak{s}}}}_{\mathcal{S}}\\
&\leq &C\|\Omega\|_{L^\rho(\mathbb{S}^{n-1})}[\omega]_{A_\delta}^\frac{1}{\mathfrak{s}\delta} \left(\mathscr{M}_\omega(|\nabla f|^{\mathfrak{s}\delta})(x)\right)^{\frac{1}{\mathfrak{s}\delta}(1-\frac{\mathfrak{q}}{\bf d})}  \|\nabla f\|_{\mathcal{M}^{\mathfrak{s}\delta, \mathfrak{q}}(\omega)}^{\frac{\mathfrak{q}}{\bf d}}.
\end{eqnarray*}
Finally, from this uniform inequality we can deduce the estimate
$$T^*_\Omega(f)(x)\leq C\|\Omega\|_{L^\rho(\mathbb{S}^{n-1})}[\omega]_{A_\delta}^\frac{1}{\mathfrak{s}\delta} \left(\mathscr{M}_\omega(|\nabla f|^{\mathfrak{s}\delta})(x)\right)^{\frac{1}{\mathfrak{s}\delta}(1-\frac{\mathfrak{q}}{\bf d})} \|\nabla f\|_{\mathcal{M}^{\mathfrak{s}\delta, \mathfrak{q}}(\omega)}^{\frac{\mathfrak{q}}{\bf d}},$$
 and this ends the proof of the Theorem \ref{Theo_Pointwise_Weighted}. \hfill $\blacksquare$
\mysection{Proof of the Corollaries}\label{Secc_ProofCorollaries}

\noindent{\bf Proof of the Corollary \ref{Corollary2}.}
\begin{itemize}
\item[1)] For the first functional inequality we use as starting point the estimate 
$$T^*_\Omega(f)(x)\leq C\|\Omega\|_{L^\rho(\mathbb{S}^{n-1})}[\omega]_{A_\delta}^\frac{1}{\mathfrak{s}\delta} \left(\mathscr{M}_\omega(|\nabla f|^{\mathfrak{s}\delta})(x)\right)^{\frac{1}{\mathfrak{s}\delta}(1-\frac{\mathfrak{q}}{\bf d})} \|\nabla f\|_{\mathcal{M}^{\mathfrak{s}\delta, \mathfrak{q}}(\omega)}^{\frac{\mathfrak{q}}{\bf d}},$$
then we apply the norm $L^r(\omega)$ to both sides of the previous expression to obtain 
$$\|T^*_\Omega(f)\|_{L^r(\omega)}\leq C_{\Omega, \omega}\left\|\left(\mathscr{M}_\omega(|\nabla f|^{\mathfrak{s}\delta})\right)^{\frac{1}{\mathfrak{s}\delta}(1-\frac{\mathfrak{q}}{\bf d})}\right\|_{L^r(\omega)} \|\nabla f\|_{\mathcal{M}^{\mathfrak{s}\delta, \mathfrak{q}}(\omega)}^{\frac{\mathfrak{q}}{\bf d}},$$
and since by the property (\ref{ProprietePuissanceLebesgue}) we have the identity 
$$\left\|\left(\mathscr{M}_\omega(|\nabla f|^{\mathfrak{s}\delta})\right)^{\frac{1}{\mathfrak{s}\delta}(1-\frac{\mathfrak{q}}{\bf d})}\right\|_{L^r(\omega)}=\left\|\mathscr{M}_\omega(|\nabla f|^{\mathfrak{s}\delta})\right\|_{L^{\frac{r}{\mathfrak{s}\delta}(1-\frac{\mathfrak{q}}{\bf d})}(\omega)}^{\frac{1}{\mathfrak{s}\delta}(1-\frac{\mathfrak{q}}{\bf d})},$$
we obtain 
$$\|T^*_\Omega(f)\|_{L^r(\omega)}\leq C_{\Omega, \omega} \left\|\mathscr{M}_\omega(|\nabla f|^{\mathfrak{s}\delta})\right\|_{L^{\frac{r}{\mathfrak{s}\delta}(1-\frac{\mathfrak{q}}{\bf d})}(\omega)}^{\frac{1}{\mathfrak{s}\delta}(1-\frac{\mathfrak{q}}{\bf d})}\|\nabla f\|_{\mathcal{M}^{\mathfrak{s}\delta, \mathfrak{q}}(\omega)}^{\frac{\mathfrak{q}}{\bf d}}.$$
Since $\frac{r}{\mathfrak{s}\delta}(1-\frac{\mathfrak{q}}{\bf d})>1$ and since the weighted maximal function $\mathscr{M}_\omega$ is bounded in the weighted Lebesgue space $L^{\frac{r}{\mathfrak{s}\delta}(1-\frac{\mathfrak{q}}{\bf d})}(\omega)$, we have
$$\|T^*_\Omega(f)\|_{L^r(\omega)}\leq C_{\Omega, \omega} \left\||\nabla f|^{\mathfrak{s}\delta}\right\|_{L^{\frac{r}{\mathfrak{s}\delta}(1-\frac{\mathfrak{q}}{\bf d})}(\omega)}^{\frac{1}{\mathfrak{s}\delta}(1-\frac{\mathfrak{q}}{\bf d})}\|\nabla f\|_{\mathcal{M}^{\mathfrak{s}\delta, \mathfrak{q}}(\omega)}^{\frac{\mathfrak{q}}{\bf d}}.$$
Applying again the property (\ref{ProprietePuissanceLebesgue}) we obtain
\begin{equation}\label{EstimationFonctionnelle1}
\|T^*_\Omega(f)\|_{L^r(\omega)}\leq C_{\Omega, \omega} \|\nabla f\|_{L^{r(1-\frac{\mathfrak{q}}{\bf d})}(\omega)}^{(1-\frac{\mathfrak{q}}{\bf d})}\|\nabla f\|_{\mathcal{M}^{\mathfrak{s}\delta, \mathfrak{q}}(\omega)}^{\frac{\mathfrak{q}}{\bf d}},
\end{equation}
but since $\frac{r}{\mathfrak{s}\delta}(1-\frac{\mathfrak{q}}{\bf d})=\mathfrak{q}$, we have 
$$\|T^*_\Omega(f)\|_{L^r(\omega)}\leq C_{\Omega, \omega} \|\nabla f\|_{L^{\mathfrak{q}}(\omega)}^{(1-\frac{\mathfrak{q}}{\bf d})}\|\nabla f\|_{\mathcal{M}^{\mathfrak{s}\delta, \mathfrak{q}}(\omega)}^{\frac{\mathfrak{q}}{\bf d}}.$$
Now, since we have the space inclusion $L^{\mathfrak{q}}(\omega)=\mathcal{M}^{\mathfrak{q}, \mathfrak{q}}(\omega)\subset \mathcal{M}^{\mathfrak{s}\delta, \mathfrak{q}}(\omega)$ (recall that $\mathfrak{s}\delta\leq \mathfrak{q}$), we can write
\begin{eqnarray*}
\|T^*_\Omega(f)\|_{L^r(\omega)}&\leq & C_{\Omega, \omega} \|\nabla f\|_{L^{\mathfrak{q}}(\omega)}^{(1-\frac{\mathfrak{q}}{\bf d})}\|\nabla f\|_{\mathcal{M}^{\mathfrak{s}\delta, \mathfrak{q}}(\omega)}^{\frac{\mathfrak{q}}{\bf d}}\leq  C_{\Omega, \omega} \|\nabla f\|_{L^{\mathfrak{q}}(\omega)}^{(1-\frac{\mathfrak{q}}{\bf d})}\|\nabla f\|_{\mathcal{M}^{\mathfrak{q}, \mathfrak{q}}(\omega)}^{\frac{\mathfrak{q}}{\bf d}}\\
&\leq & C_{\Omega, \omega} \|\nabla f\|_{L^{\mathfrak{q}}(\omega)}^{(1-\frac{\mathfrak{q}}{\bf d})}\|\nabla f\|_{L^{\mathfrak{q}}(\omega)}^{\frac{\mathfrak{q}}{\bf d}}\leq C_{\Omega, \omega} \|\nabla f\|_{L^{\mathfrak{q}}(\omega)},
\end{eqnarray*}
and the first functional inequality is thus proven.

\item[2)] For the second functional inequality, we follow the same steps used to proof the estimate (\ref{EstimationFonctionnelle1}) and we have
$$\|T^*_\Omega(f)\|_{L^r(\omega)}\leq C_{\Omega, \omega} \|\nabla f\|_{L^{r(1-\frac{\mathfrak{q}}{\bf d})}(\omega)}^{(1-\frac{\mathfrak{q}}{\bf d})}\|\nabla f\|_{\mathcal{M}^{\mathfrak{s}\delta, \mathfrak{q}}(\omega)}^{\frac{\mathfrak{q}}{\bf d}},$$
now, since $r(1-\frac{\mathfrak{q}}{\bf d})=\sigma$, we obtain
$$\|T^*_\Omega(f)\|_{L^r(\omega)}\leq C_{\Omega, \omega} \|\nabla f\|_{L^{r(1-\frac{\mathfrak{q}}{\bf d})}(\omega)}^{(1-\frac{\mathfrak{q}}{\bf d})}\|\nabla f\|_{\mathcal{M}^{\mathfrak{s}\delta, \mathfrak{q}}(\omega)}^{\frac{\mathfrak{q}}{\bf d}}\leq C_{\Omega, \omega} \|\nabla f\|_{L^{\sigma}(\omega)}^{(1-\frac{\mathfrak{q}}{\bf d})}\|\nabla f\|_{\mathcal{M}^{\mathfrak{s}\delta, \mathfrak{q}}(\omega)}^{\frac{\mathfrak{q}}{\bf d}},$$
which proves the second inequality of the Corollary \ref{Corollary2}. 

\item[3)]  The proof of the last inequality of the Corollary \ref{Corollary2} is very similar. Indeed, we write
$$T^*_\Omega(f)(x)\leq C\|\Omega\|_{L^\rho(\mathbb{S}^{n-1})}[\omega]_{A_\delta}^\frac{1}{\mathfrak{s}\delta} \left(\mathscr{M}_\omega(|\nabla f|^{\mathfrak{s}\delta})(x)\right)^{\frac{1}{\mathfrak{s}\delta}(1-\frac{\mathfrak{q}}{\bf d})} \|\nabla f\|_{\mathcal{M}^{\mathfrak{s}\delta, \mathfrak{q}}(\omega)}^{\frac{\mathfrak{q}}{\bf d}},$$
we apply now the weigthed Morrey norm $\mathcal{M}^{\mathfrak{a}, \mathfrak{b}}(\omega)$ to both sides if the previous pointwise inequality to obtain
\begin{eqnarray*}
\|T^*_\Omega(f)\|_{\mathcal{M}^{\mathfrak{a}, \mathfrak{b}}(\omega)}&\leq &C\|\Omega\|_{L^\rho(\mathbb{S}^{n-1})}[\omega]_{A_\delta}^\frac{1}{\mathfrak{s}\delta} \left\|\left(\mathscr{M}_\omega(|\nabla f|^{\mathfrak{s}\delta})\right)^{\frac{1}{\mathfrak{s}\delta}(1-\frac{\mathfrak{q}}{\bf d})}\right\|_{\mathcal{M}^{\mathfrak{a}, \mathfrak{b}}(\omega)} \|\nabla f\|_{\mathcal{M}^{\mathfrak{s}\delta, \mathfrak{q}}(\omega)}^{\frac{\mathfrak{q}}{\bf d}}\\
&\leq &C_{\Omega, \omega} \left\|\left(\mathscr{M}_\omega(|\nabla f|^{\mathfrak{s}\delta})\right)\right\|_{\mathcal{M}^{\mathfrak{a}({\frac{1}{\mathfrak{s}\delta}(1-\frac{\mathfrak{q}}{\bf d})}), \mathfrak{b}({\frac{1}{\mathfrak{s}\delta}(1-\frac{\mathfrak{q}}{\bf d})})}(\omega)}^{\frac{1}{\mathfrak{s}\delta}(1-\frac{\mathfrak{q}}{\bf d})} \|\nabla f\|_{\mathcal{M}^{\mathfrak{s}\delta, \mathfrak{q}}(\omega)}^{\frac{\mathfrak{q}}{\bf d}},
\end{eqnarray*}
where we applied the property (\ref{ProprietePuissanceMorrey}). Recalling that $1< \mathfrak{a}({\frac{1}{\mathfrak{s}\delta}(1-\frac{\mathfrak{q}}{\bf d})})\leq  \mathfrak{b}({\frac{1}{\mathfrak{s}\delta}(1-\frac{\mathfrak{q}}{\bf d})})<+\infty$, then the weighted maximal function $\mathscr{M}_\omega$ is bounded in this weighted Morrey space and we can write
\begin{eqnarray*}
\|T^*_\Omega(f)\|_{\mathcal{M}^{\mathfrak{a}, \mathfrak{b}}(\omega)}&\leq &C_{\Omega, \omega} \left\||\nabla f|^{\mathfrak{s}\delta}\right\|_{\mathcal{M}^{\mathfrak{a}({\frac{1}{\mathfrak{s}\delta}(1-\frac{\mathfrak{q}}{\bf d})}), \mathfrak{b}({\frac{1}{\mathfrak{s}\delta}(1-\frac{\mathfrak{q}}{\bf d})})}(\omega)}^{\frac{1}{\mathfrak{s}\delta}(1-\frac{\mathfrak{q}}{\bf d})} \|\nabla f\|_{\mathcal{M}^{\mathfrak{s}\delta, \mathfrak{q}}(\omega)}^{\frac{\mathfrak{q}}{\bf d}}\\
&\leq & C_{\Omega, \omega} \|\nabla f\|_{\mathcal{M}^{\mathfrak{a}(1-\frac{\mathfrak{q}}{\bf d}), \mathfrak{b}(1-\frac{\mathfrak{q}}{\bf d})}(\omega)}^{(1-\frac{\mathfrak{q}}{\bf d})} \|\nabla f\|_{\mathcal{M}^{\mathfrak{s}\delta, \mathfrak{q}}(\omega)}^{\frac{\mathfrak{q}}{\bf d}},
\end{eqnarray*}
and this proves the last inequality of the Corollary \ref{Corollary2}. \hfill$\blacksquare$\\
\end{itemize}
\noindent {\bf Acknowledgment.} This work was supported by the GDRI ECO-Math.\\

\noindent {\bf Conflict of interest.} We declare that we do not have any commercial or associative interest that represents a conflict of interest in connection with the work submitted.


\end{document}